\documentclass[11pt]{article}

\usepackage{amsmath, amssymb, amsthm}

\usepackage[utf8]{inputenc}

\usepackage{authblk}

\usepackage{geometry}
\usepackage{microtype}

\usepackage[colorlinks=true, linkcolor=blue, citecolor=blue, urlcolor=blue]{hyperref}

\usepackage{enumitem}

\usepackage{tikz}
\usetikzlibrary{shapes.geometric, patterns, arrows.meta, decorations.pathreplacing, calc}
\tikzset{>=Stealth} 
\definecolor{lightorange}{rgb}{1, 0.85, 0.7}
\definecolor{lightgray}{rgb}{0.9, 0.9, 0.9}

\numberwithin{equation}{section}

\newtheorem{corollary}{Corollary}
\newtheorem{theorem}{Theorem}

\newcommand{\esssup}{\operatorname*{ess\,sup}}

\title{\textbf{Dynamic Refinement of Pressure Decomposition in Navier-Stokes Equations}}
\author[1]{P. G. Fernández-Dalgo}
\affil[1]{Basque Center for Applied Mathematics, BCAM, Bilbao, España \\ \textit{Email address:} \texttt{pfernandez@bcamath.org}}
\date{} 

\begin{document}

\maketitle

\begin{abstract}
In this work, the local decomposition of pressure in the Navier-Stokes equations is dynamically refined to prove that a relevant critical energy of a suitable Leray-type solution inside a backward paraboloid---regardless of its aperture---is controlled near the vertex by a critical behavior confined to a neighborhood of the paraboloid's boundary. This neighborhood excludes the interior near the vertex and remains separated from the temporal profile of the vertex, except at the vertex itself. Then, we present a refined scaling invariant regularity result. 
\end{abstract}



\section{Introduction}

The recent preprint by Coiculescu and Palasek \cite{CoiculescuPalasek2025} presents a remarkable result showing non-uniqueness of smooth solutions to the incompressible Navier--Stokes equations with an explicit multi-scale
initial datum which belongs to the critical space \( BMO^{-1} \). This breakthrough contrasts sharply with the still open question of regularity in large critical spaces.

Continuing the investigation conducted in \cite{BarkerFernandezPrange2024} to address the endpoint critical case for the regularity problem of the three dimensional Navier-Stokes equations,
\begin{equation}
\label{NS}
    \partial_t v - \Delta v + v \cdot \nabla v + \nabla p = 0, \quad \nabla \cdot v = 0,
\end{equation}
we aim to relax the following hypothesis, written in terms of the function $\theta_a(s) = \sqrt{-as}$,
\begin{equation*}
    \operatorname*{\,limsup}_{s \to 0^+} \|v(\cdot, s)\|_{L^{3}(B(\sqrt{a}) \setminus B(\theta_{a}(s)))} < + \infty,
\end{equation*}
which was imposed in \cite[Theorem A]{BarkerFernandezPrange2024} on the exterior of the paraboloid
\begin{align*}
\label{Pa}
P_{\sqrt{a}}(-1)  &:= \, \displaystyle \bigcup\limits_{s \in \left(-1, 0\right)}  \{x :  |x|=\theta_{a}(s) \} \, \times \, \{ \tau \},
\end{align*}
with a small aperture $\sqrt{a} \approx 1$, to propagate regularity up to $(0,0)$.
To achieve this, we will make a dynamic partition of spatial scales and rely on $L^\infty_t L^{3,\infty}$ behavior in a neighborhood of the paraboloid's boundary with aperture $\sqrt{a}N$ with $N \gg 1$, this neighborhood remains separated from the temporal profile.  
\noindent
Before introducing the context and presenting our contributions, let us clarify the notations used throughout this text.

\subsection*{Notations.}
We denote by \( B(r) \) the ball centered at the origin with radius \( r > 0 \). 
Throughout this work, \( C \) will represent a positive universal constant, which may vary between different expressions or lines. Importantly, \( C \) remains independent of the parameters \( a \), \( \gamma \), or \( N \). 
When a constant depends on specific parameters \( b_1, \ldots, b_k \), we denote it by \( C_{b_1, \ldots, b_k} \). The constants that need to be referenced in subsequent calculations will be distinguished by a superscript, such as \( C^{p_2} \).

\begin{figure}[htbp]
\centering
\begin{tikzpicture}[scale=3]


\draw[thick, domain=-1:1, samples=100] 
    plot(\x, {-\x*\x});

\draw[thick, orange, domain=-0.25:0.25, samples=100] 
    plot(4*\x, {-\x*\x});

\draw[->, thick] (-1.2,0) -- (1.2,0) node[right] {$x$};
\draw[->, thick] (0,-1.1) -- (0,0.2) node[above] {$t$};
\draw[-] (-1,-1.1) -- (-1,0.1) node[above] {$x=\sqrt{a}$};

\node[orange] at (2.2, -0.1) {\small Parabolic scales of order $\sqrt{a}N$: $|x| = \sqrt{-at}N$, $N>>1$};
\node[black] at (1.6, -0.5) {\small Parabolic scales of order $\sqrt{a}$: $|x| = \sqrt{-at}$, $\sqrt{a} \approx 1$};

\draw[fill=blue] (0,0) circle [radius=0.02];

\end{tikzpicture}
\caption{Parabolic scales}
\end{figure}
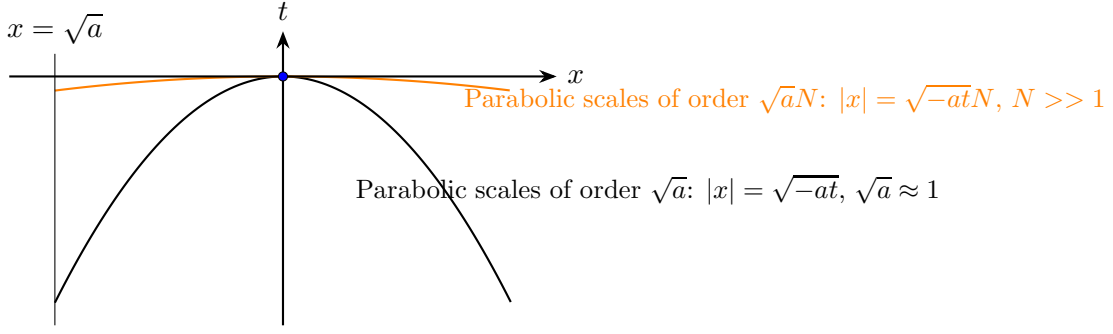

\section*{}
Key ingredients for propagating regularity from the exterior of a paraboloid with small aperture in the endpoint case, as discussed in \cite[Theorem A]{BarkerFernandezPrange2024}, include the introduction of a second aperture parameter \( N \) into the critical energy balance \eqref{balance} and, along a wider paraboloid, controlling the \footnote{Criticality with respect to the scaling invariance for the Navier-Stokes equations is explained in Appendix \ref{critical}.}critical quantities \( f \) and \( g(\cdot,0) \), given by:

\begin{equation}
\label{fg}
    f(s) := \frac{1}{\theta_a(s)} \|\Psi_N v(\cdot, s)\|_{L^2}^2, \quad g_\gamma(s, t) := \int_s^t \frac{\theta_a(\tau)^{\gamma-1}}{\theta_a(s)^{\gamma+1}} f(\tau) \, d\tau, \quad \gamma >-1,
\end{equation}
where we consider \( \Psi_N(x,t) := \varphi_N\left(\frac{x}{\theta_a(t)}\right) \), and \( \varphi_N \in C^\infty_c(\mathbb{R}^3) \) is a smooth, compactly supported, radially decreasing and positive function fulfilling \( \varphi_N(x) = 1 \) within the ball \( B(N) \), its support satisfying \( \operatorname{supp} \varphi_N \subset B(N + 1) \), and its gradient being bounded \( \|\nabla \varphi_N\|_{L^\infty} \lesssim 1 \) uniformly on N. Consequently, \( \Psi_N(x,t) \) inherits these properties, which allows us to get for $\tau<0$,
\begin{equation*}
    \partial_t \left( \Psi_N ^2 (x,\cdot) \right) (\tau) \leq 0,
\end{equation*}
\begin{align*}
    \text{supp} (\Psi_N (\cdot,\tau)) \subset B \left( N \theta_a(\tau) \right), \text{ and }
\end{align*}
\begin{align}
\label{suppgradpsi}
    \text{supp} (\nabla \Psi_N (\cdot,\tau)) &\subset B \left( N \theta_a(t) \right) \setminus B \left( (N+1) \theta_a(\tau) \right) .
\end{align}
Thus, the integration domain of \( g_\gamma(s) := g_\gamma(s,0) \) is contained within the region inside the paraboloid (with aperture \( \sqrt{a}(N+1) \)),
\begin{align*}
\label{domintg}
    \displaystyle \bigcup\limits_{\tau \in \left(-s, 0\right)}  B \left( (N+1) \theta_a(\tau) \right) \, \times \, \{ \tau \} .
\end{align*}

\noindent
Reinjecting the bounds for \( f \) and \( g_\gamma \), they also obtain the bound for \( h_\gamma(s) := h_\gamma(s,0) \), defined by:
\begin{equation}
\label{h}
    h_\gamma(s, t) := \int_s^t \frac{\theta_a(\tau)^{\gamma}}{\theta_a(s)^{\gamma+1}} \| \nabla (\Psi_N v(\cdot, \tau )) \|_{L^2}^2 \, d\tau, \quad \gamma >-1.
\end{equation}

\noindent
To bound these critical quantities through the energy balance, \cite{BarkerFernandezPrange2024} performs a spatially localized decomposition of the pressure. The authors observe that applying a parabolic cut-off for the spatial localization very close to the aperture \( \sqrt{a} \) (and hence far from the aperture \( \sqrt{a}N \)) allows them to absorb the non-local effect of the pressure and bound the energy using a Gronwall type inequality, following the ideas in \cite{Neustupa2014}. For this reason, the cut-off is applied very close to the aperture \( \sqrt{a} \).

\noindent
We will prove that, to bound the key term in the balance
\begin{equation}
\label{prepressure}
     \int_s^t \frac{\theta_a(\tau)^{\gamma - 1}}{\theta_a(s)^{\gamma + 1}} \int_{\text{supp} (\nabla \Psi_N (\cdot,\tau))} |p v \nabla \Psi_N| \, dx \, d\tau,
\end{equation}
it is more effective to employ a dynamic localization and consider a parabolic cut-off for the localization that optimally absorbs the non-local effect of the pressure. 

Let us now specify our modification for this technical point (with respect to what was done in \cite{BarkerFernandezPrange2024}). The key lies in introducing a dynamic localization adapted to the term \eqref{prepressure} and to the context, for instance the context of \footnote{We refer to \cite[Definition 6.2]{FernandezDalgo2021}. Let us mention that we use only the fact that $v$ is locally $L^\infty L^2$, $\nabla v$ is locally $L^2 L^2$ and $p$ locally $L^{\frac{3}{2}}L^{\frac{3}{2}}$, in particular, to control the $p_2$ part of $p$ in Section \ref{rhsbalance}, we use $v$ is locally $L^3L^3$.}suitable Leray-type solutions is considered in Theorem \ref{premain} Part A). We fix the dynamic localization \footnote{For the Part B) in Theorem \ref{premain} we will consider a parabolic localization} as follows: 
\begin{equation*}
    \eta_N(x,\tau) = \varphi_{N+2} \left( \frac{x}{\theta_a(\tau)^{1/2}} \right) .
\end{equation*}
Then, we have 
\begin{align}
\label{suppgradeta}
    \text{supp} (\nabla \eta_N(\cdot, \tau)) &\subset B \left( (N+3) \theta_a(\tau)^{1/2} \right) \setminus B \left( (N+2) \theta_a(\tau)^{1/2} \right)
\end{align}
and the following bounds hold: \( |\nabla \eta_N| \lesssim \frac{C}{\theta(t)^{1/2}} \) and \( |\nabla^2 \eta_N| \lesssim \frac{C}{\theta(t)} \). 

\begin{figure}[htbp]
\centering
\begin{tikzpicture}[scale=2]
    \fill[lightorange] 
        plot[domain=-2.65:2.65, samples=100] (\x, {- (\x*\x*\x*\x)/(9*9) }) -- 
        plot[domain=2.3:-2.3, samples=100] (\x, {- (\x*\x) / (3*3)}) -- 
        cycle;
    \fill[lightgray] 
        plot[domain=-2.5:2.5, samples=100] (\x, {- 0.6 }) -- 
        plot[domain=2.3:-2.3, samples=100] (\x, {- (\x*\x) / (3*3)}) -- 
        cycle;
        
    \draw[->, thick] (0,-1) -- (0,0.3) node[above] {$t$};
    
    \draw[thick, orange] plot[domain=-2.3:2.3, samples=100] (\x, {- (\x*\x) / (3*3)}) node[below] {};
    \draw[thick, red] plot[domain=-2.65:2.65, samples=100] (\x, {- (\x*\x*\x*\x)/(9*9)}) node[right] {};
    
    \node[red] at (3.6, -0.1) {$|x|=(N+3)(-at)^{1/4} = (N+3) \theta_a(t)^{\frac{1}{2}}$};
    
    \node[orange] at (3.1, -0.8) {$|x|= \frac{N}{2}\sqrt{-at} = \frac{N}{2} \theta_a(t)$};

\end{tikzpicture}
\caption{Regions in dynamic decomposition with $a=1$ and $N=6$}
\end{figure}
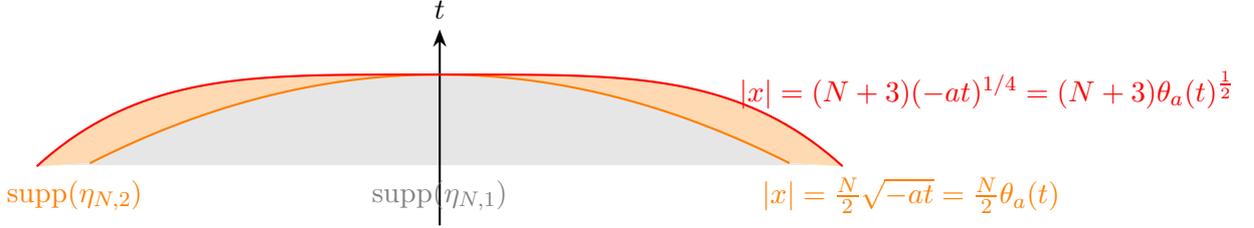

\noindent
Then, using the identity
\[
\eta_N(x, \tau)p(x, \tau) = -\frac{1}{4\pi} \int_{\mathbb{R}^3} \frac{1}{|x-y|} [\Delta (\eta_N p)](y, \tau) \, dy,
\]
the fact that \( \Delta p = -\partial_i \partial_j (v_i v_j) \), and some standard computations one can verify the pressure decomposition
\begin{equation}
\label{decomposition}
    \eta_N(x)p(x, \tau) := p_1(x, \tau) + p_2(x, \tau) + p_3(x, \tau),
\end{equation}
where
\[
p_1(x, \tau) := \frac{1}{4\pi} \int_{B((\frac{N}{2}+1))\theta_a(\tau))} \frac{\partial^2}{\partial y_i \partial y_j} \left( \frac{1}{|x-y|} \right) [\eta_N v_i v_j](y, \tau) \, dy,
\]

\[
p_2(x, \tau) := \frac{1}{4\pi} \int_{B\left((N+3) \theta_a(\tau)^{1/2}\right) \setminus B((\frac{N}{2}+1))\theta_a(\tau))} \frac{\partial^2}{\partial y_i \partial y_j} \left( \frac{1}{|x-y|} \right) [\eta_N v_i v_j](y, \tau) \, dy
\]
and
\begin{align*}
    p_3(x, \tau) :=& \frac{1}{2\pi} \int_{\text{supp} (\nabla \eta_N(\cdot, \tau))} \frac{x_i - y_i}{|x-y|^3} \left( \frac{\partial \eta_N}{\partial y_j} v_i v_j \right)(y, \tau) \, dy \\
    & + \frac{1}{4\pi} \int_{\text{supp} (\nabla \eta_N(\cdot, \tau))} \frac{1}{|x-y|} \left( \frac{\partial^2 \eta_N}{\partial y_i \partial y_j} v_i v_j \right)(y, \tau) \, dy \\
    & + \frac{1}{2\pi} \int_{\text{supp} (\nabla \eta_N(\cdot, \tau))} \frac{x_i - y_i}{|x-y|^3} \left( \frac{\partial \eta_N}{\partial y_j} p \right)(y, \tau) \, dy \\
    &+ \frac{1}{4\pi} \int_{\text{supp} (\nabla \eta_N(\cdot, \tau))} \frac{1}{|x-y|} (\Delta \eta_N p)(y, \tau) \, dy.
\end{align*}
In the part \( p_1 \), from \eqref{suppgradpsi} we observe that
\begin{equation}
    \text{dist}\big(\text{supp}(\nabla \Psi_N(\cdot,\tau) ), B((\frac{N}{2}+1) \theta_a(\tau)) \big) \geq  (\frac{N}{2}-1)\theta_a(\tau),
\end{equation}
which immediately implies that for 
$  x \in \text{supp}(\nabla \Psi_N(\cdot,\tau))$,
\begin{equation}
\label{kernelp1}
    |p_1(x, \tau)| \leq \frac{C}{N^3 \theta_a(\tau)^3} \int_{B(\frac{N}{2}\theta_a(\tau))} |v(y, \tau)|^2 \, dy .
\end{equation}
For \( p_3 \), we observe that for \( \tau \in \left( -\frac{1}{4a} , 0 \right) \) (and hence \( \theta_a(\tau) < \frac{1}{2} \)), 
\begin{align*}
    \text{dist} & \big( \text{supp}(\nabla \Psi_N(\cdot,\tau)), \text{supp} (\nabla \eta_N(\cdot, \tau)) \big) \nonumber\\
    &\geq (N+2) \theta_a(\tau)^{1/2} - (N+1) \theta_a(\tau) \nonumber\\
    &\geq C \theta_a(\tau)^{1/2}.
\end{align*}
which implies that
\begin{equation}
\label{kernelp3}
    |p_3(x, \tau)| \leq \frac{C}{\theta(\tau)^\frac{3}{2}} \int_{\text{supp} (\nabla \eta_N(\cdot, \tau))} \big(|v|^2 + |p| \big)(y, \tau) \, dy.
\end{equation}
For \( p_2 \), the singularity arises, and we invoke Calderón-Zygmund's theorem, which provides,
\begin{equation}
\label{kernelp2}
    \|p_2(\cdot, \tau)\|_{L^{\frac{3}{2},\infty}(\text{supp}(\nabla \Psi_N(\cdot,\tau)))} \leq C \| v(\cdot, \tau)\|^2_{L^{3,\infty}(B((N+3)\theta_a(\tau)^{1/2}) \setminus B((\frac{N}{2}+1)\theta_a(\tau)))}.
\end{equation}

\noindent
We will verify that, with this decomposition, the rest of the regularity analysis proceeds correctly essentially from the hypothesis
\begin{equation}
\label{hip2main}
    \operatorname*{ess\,sup}_{s \in (-c,0)} \|v(\cdot, s)\|_{L^{3, \infty}(B(b) \,\, \cap \,\, B((N+3) \theta_a(s)^{1/2}) \setminus B(\frac{N}{2} \theta_{a}(s)))} \leq + \infty.
\end{equation}

Moreover, we will explicitly provide the details to consider the constant \( a \) within the complete range \( a \in \left(0, 4 \lambda_S(B(1))\right) \), where \( \lambda_S(B(1)) > \pi^2 \) is the first eigenvalue of the Dirichlet-Stokes operator on \( B(1) \), as initially considered by \cite{Neustupa2014}.

This restriction for $a$ arises from the need to gain critical \( L^2 \) information at the inner scales (within the paraboloid of aperture \( \sqrt{a}N \)) through critical \( L^2 \) information on the gradient of the fluid at those scales. The use of a Poincaré inequality to achieve this imposes a bound on the aperture parameter \( a \). Consequently, we are compelled to assume information at intermediate scales (between \( \sqrt{a} \) and \( \sqrt{a}N \)). This critical information is \( L^2 \) in nature for the velocity, and thus it is significantly weaker than the information assumed in \cite[Theorem A]{BarkerFernandezPrange2024}.

\noindent
Thus, our main contribution is the following:

\begin{theorem}
\label{premain}
Let \footnote{Our hypotheses in A) are not scaling-invariant due to \eqref{hip2main}, as the set \( B((N+3) \theta_a(s)^{1/2}) \setminus B(\frac{N_0}{2} \theta_{a}(s)) \) is not parabolic. Thus, to address the general case, we introduce the parameters \( a, b, c \), where \( a \) determines the aperture of the paraboloid, \( b \) specifies the spatial localization of the hypotheses, and \( c \) defines the regularity interval of the solution \( v \).} \( a, b, c > 0 \) with $a \in \left(0, 4 \lambda_S(B(1))\right) $ and let \( M \geq 1 \). Let \( v \) be a suitable Leray-type solution\footnote{see \cite[Definition 6.2]{FernandezDalgo2021}} to the Navier-Stokes equations \eqref{NS} in \( \mathbb{R}^3 \times (-c, 0) \), such that \( v \in C^\infty((-c, T); C^\infty(B(b))) \) for all \( T \in (-c, 0) \).

\noindent
Then, 
\begin{enumerate}[label=\Alph*)]
        
    \item there exists \( N_0 > 0 \) such that for all \( N \geq N_0 \), if 
        \begin{equation}
        \label{hip1main}
            \operatorname*{ess\,sup}_{s \in (-c,0)}  \frac{1}{N \theta_a(s)} \int_{ B(b) \,\, \cap \,\, B((N+1) \theta_{a}(s)) \setminus B(\theta_a(s))} |v(x, s)|^2 \, dy  \leq CM,
        \end{equation}
        and
        \begin{equation}
        \label{hip2main}
            \operatorname*{ess\,sup}_{s \in (-c,0)} \|v(\cdot, s)\|_{L^{3, \infty}(B(b) \,\, \cap \,\, B((N+3) \theta_a(s)^{1/2}) \setminus B(\frac{N}{2} \theta_{a}(s)))} \leq M
        \end{equation}
        then the ($N$-dependent) scaling invariant function \( f \) defined in \eqref{fg} is bounded. Moreover, for \( \gamma > -1 \), the ($N$-dependent) scaling invariant functions \( g_\gamma = g_\gamma(\cdot,0) \) and \( h_\gamma = h_\gamma(\cdot,0) \) (see \eqref{h}) are also bounded.

        \item there exists \( N_0 > 0 \) such that for all \( N \geq N_0 \), if 
        \begin{equation}
        \label{hip1maincri}
            \operatorname*{ess\,sup}_{s \in (-c,0)}  \frac{1}{N \theta_a(s)} \int_{ B(b) \,\, \cap \,\, B((2N+2) \theta_{a}(s)) \setminus B(\theta_a(s))} |v(x, s)|^2 \, dy  \leq CM,
        \end{equation}
        \begin{equation}
        \label{hip2maincri}
            \operatorname*{ess\,sup}_{s \in (-c,0)} \|v(\cdot, s)\|_{L^{3, \infty}(B(b) \,\, \cap \,\, B((2N+2) \theta_a(s)) \setminus B(\frac{N}{2} \theta_{a}(s)))} \leq M
        \end{equation}
        and
        \begin{equation}
        \label{weakcrip}
            \operatorname*{ess\,sup}_{s \in (-c,0)}  \frac{1}{N \theta_a(s)} \int_{ B(b) \,\, \cap \,\, B((2N+2) \theta_{a}(s)) \setminus B((2N+1)\theta_a(s))} |p(x, s)| \, dy  \leq CM,
        \end{equation}
        then the functions \( f \), \( g_\gamma \) and \( h_\gamma \) are bounded, for \( \gamma > -1 \), and more precisely
        \begin{equation*}
        \label{fcri}
            \frac{1}{N \theta_a(s)} \int_{B(\frac{b}{16})} | \Psi_N v(x, s)|^2 dx \leq C_{a, \gamma, M},
        \end{equation*}
        \begin{equation*}
        \label{gcri}
            \frac{1}{N} \int_s^t \frac{\theta_a(\tau)^{\gamma-2}}{\theta_a(s)^{\gamma+1}} \int_{B(\frac{b}{16})} | \Psi_N v(x, \tau )|^2 dx \, d\tau \leq C_{a, \gamma, M}
        \end{equation*}
        and
        \begin{equation*}
        \label{hcri}
            \frac{1}{N} \int_s^t \frac{\theta_a(\tau)^{\gamma}}{\theta_a(s)^{\gamma+1}} \int_{B(\frac{b}{16})} | \nabla (\Psi_N v(\cdot, \tau )) |^2 (x) dx \, d\tau \leq C_{a, \gamma, M},
        \end{equation*}
        where $C_{a, \gamma, M}$ does not depend on $N$.

\end{enumerate}
\end{theorem}

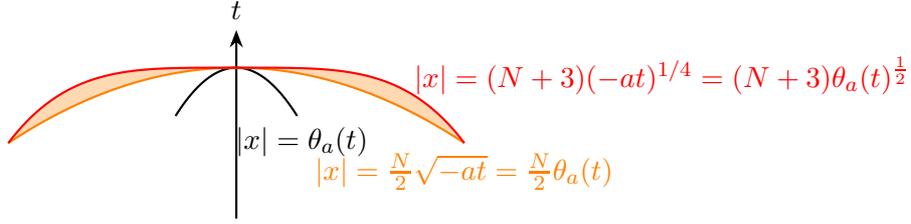
\begin{figure}[htbp]
\centering
\begin{tikzpicture}[scale=1]
    \fill[lightorange] 
        plot[domain=-3:3, samples=100] (\x, {- (\x*\x*\x*\x)/(9*9) }) -- 
        plot[domain=3:-3, samples=100] (\x, {- (\x*\x) / (3*3)}) -- 
        cycle;
        
    \draw[->, thick] (0,-2) -- (0,0.5) node[above] {$t$};
    
    \draw[thick] plot[domain=-0.8:0.8, samples=100] (\x, {-1*\x*\x}) node[below] {$\,\, |x|= \theta_a(t)$};
    \draw[thick, orange] plot[domain=-3:3, samples=100] (\x, {- (\x*\x) / (3*3)}) node[below] {$|x|= \frac{N}{2}\sqrt{-at} = \frac{N}{2} \theta_a(t)$};
    \draw[thick, red] plot[domain=-3:3, samples=100] (\x, {- (\x*\x*\x*\x)/(9*9)}) node[right] {};
    
    \node[red] at (5.6, -0.1) {$|x|=(N+3)(-at)^{1/4} = (N+3) \theta_a(t)^{\frac{1}{2}}$};
\end{tikzpicture}
\caption{Non scaling invariant region in Theorem \ref{premain} taking $a=1$ and $N=6$}
\end{figure}

\begin{corollary}
\label{inter}
Let \( a, b, c > 0 \). Let \( v \) be a suitable Leray-type solution to the Navier-Stokes equations \eqref{NS} in \( \mathbb{R}^3 \times (-c, 0) \), such that \( v \in C^\infty((-c, T); C^\infty(B(b))) \) for all \( T \in (-c, 0) \).
For $N>0$ let us assume the ($N$-dependent) function $f$ is bounded. Then, for \( \gamma > -1 \), the ($N$-dependent) functions \( g_\gamma = g_\gamma(\cdot,0) \) and \( h_\gamma = h_\gamma(\cdot,0) \) are bounded. Moreover,  
\begin{equation*}
    \esssup_{s \in \left(-c, 0\right)} \frac{1}{\theta_a(s)^2} \int_s^0 \int_{B(b) \, \cap \,  B(N \theta_a(s))} |v(x, \tau)|^3 dx d\tau < C_N,
\end{equation*}
and thus we can extract a sequence of times $t_k \uparrow 0$ such that $\int_{B(N \theta_a(t_k))} | v(\cdot, t_k) |^3 \leq C_N.$ 
\end{corollary}

Our strong localized version of the critical energy bound inside a paraboloid allows us to establish a more general statement for \cite[Theorem A]{BarkerFernandezPrange2024}.

\begin{corollary}
\label{main}
Assume the parameters $a,b,c,M,N_0,N$ and the solution $v$ of \eqref{NS} satisfy the same conditions as in Theorem \ref{premain}, part A) or B). Then, if $v$ fulfills the following condition:
    \begin{equation}
    \label{hipgen}
        \esssup_{s \in \left(-c, 0\right)} \frac{1}{\theta_a(s)^2} \int_s^0 \int_{B(b) \setminus B(N \theta_a(s))} |v(x, \tau)|^3 dx d\tau < + \infty,
    \end{equation}
we conclude $(0,0)$ is a regular point.
\end{corollary}

\subsection{Outline of the paper.}
Section \ref{energy} recalls the energy balance in the interior of paraboloids and local estimates for the pressure terms (with the new dynamic decomposition). Section \ref{boundedness} is dedicated to proving Theorem \ref{premain}, while Section \ref{conclusions} concludes the proof of Corollary \ref{inter} and \ref{main}.

\section{Local energy in the interior of paraboloids}
\label{energy}
Let us briefly recall the method for obtaining the energy balance on the paraboloid. Consider \( N > 0 \) and \( s \) satisfying
\begin{equation}
\label{I}
    s \in  (-\frac{c}{2},0) \cap (-\frac{1}{2a},0) \cap \left(-\frac{b^4}{a 4^4 (N+3)^4}, 0\right) =: I_{a,b,c,N},
\end{equation}
so that for $\tau \in (s,0)$,
\begin{equation*}
    B((N+1)\theta_a(\tau)) \subset B((N+3)\theta_a(\tau)^{1/2}) \subset B(\frac{b}{4}) .
\end{equation*}
We test the Navier-Stokes equations \eqref{NS} with \footnote{$\Psi_N$ is defined in \eqref{fg} together with $f$ and $g_\gamma$.}
\[
    \theta_a^\gamma(\tau) \Psi_N^2(x,\tau), \quad \tau < 0,
\]
and  we integrate over \( \mathbb{R}^3 \times (s, t) \), where \( t \in (s, 0) \). Following \cite{BarkerFernandezPrange2024}, this leads to:
\begin{align}
\label{balance}
    \frac{\theta_a(t)^\gamma}{\theta_a(s)^{\gamma + 1}}& \|\Psi_N v(\cdot, t)\|_{L^2}^2 + \frac{a\gamma}{2} \int_s^t \frac{\theta_a(\tau)^{\gamma - 2}}{\theta_a(s)^{\gamma + 1}} \|\Psi_N v\|_{L^2}^2 \, d\tau + 2 \int_s^t \frac{\theta_a(\tau)^\gamma}{\theta_a(s)^{\gamma + 1}} \|\nabla (\Psi_N v)\|_{L^2}^2 \, d\tau \nonumber \\
    &- \int_s^t \frac{\theta_a(\tau)^{\gamma - 2}}{\theta_a(s)^{\gamma + 1}} \int_{\mathbb{R}^3} \left( \frac{1}{2} ax \cdot \nabla \Psi_N^2 \right) |v|^2 \, dx \, d\tau \nonumber \\
    = \frac{1}{\theta_a(s)} & \|\Psi_N v(\cdot, s)\|_{L^2}^2 + \int_s^t \frac{\theta_a(\tau)^\gamma}{\theta_a(s)^{\gamma + 1}} \int_{\mathbb{R}^3} \big[ 2|\nabla \Psi_N|^2 |v|^2 + (|v|^2 + 2p)(v \cdot \nabla \Psi_N^2) \big] \, dx \, d\tau .
\end{align}
Without extra assumptions on \footnote{If we assume \( f \) to be bounded (which will be the case later in the proof of Theorem \ref{main}), the terms on the right-hand side of \eqref{balance} are well-defined for \( \gamma > -1 \), as we will see in Section \ref{rhsbalance}.} \( f \), we need to assume \( \gamma > 0 \) to give sense to the energy terms on the left hand side.
\noindent
In addition to the observations in \cite{BarkerFernandezPrange2024}, we note that the term
\begin{equation}
\label{posi}
    - \int_s^t \frac{\theta_a(\tau)^{\gamma - 2}}{\theta_a(s)^{\gamma + 1}} \int_{\mathbb{R}^3} \left( \frac{1}{2} ax \cdot \nabla \Psi_N^2 \right) |v|^2 \, dx \, d\tau,
\end{equation}
placed on the left-hand side is positive because
\begin{align*}
    0 &\geq \partial_t \left( \Psi_N ^2 (x,\cdot) \right) (\tau) = 2 \Psi_N (x,\tau) \, \nabla \varphi_N \left( \frac{x}{\theta_a(\tau)} \right) \, \, \cdot \left( \frac{a}{2 \, \theta_a(\tau)^3} x \,   \right) = \frac{a}{2} \theta_a(\tau)^{-2}  2 \Psi_N  \nabla \Psi_N \cdot x \\
    &= \theta_a(\tau)^{-2} \frac{1}{2} ax \cdot   \nabla \Psi_N^2. 
\end{align*}
Thus, we do not need to control \eqref{posi}. While this term is not particularly challenging to handle, we prefer to focus our study on the key terms.

Note that \eqref{balance} can be expressed in terms of \( f \) and \( g_\gamma \), as defined in \eqref{fg}. Moreover, the more subtle quantity \( f \) can be expressed in terms of \( g_\gamma \) in the following way:
\begin{equation}
\label{fgeq}
f(s) =  a \frac{(\gamma + 1)}{2} g_\gamma(s, t) - \theta_a(s)^{-2} \frac{\partial g_\gamma}{\partial s}(s, t) .
\end{equation}
To verify this, we simply calculate the \( s \)-derivative of \( g_\gamma(\cdot, t) \),
\begin{align*}
    \frac{\partial g_\gamma}{\partial s}(s, t) &= a (\gamma + 1) \theta_a(s)^{-(\gamma+3)} \int_s^t \theta_a(\tau)^{\gamma-2} \|\Psi_N v(\cdot, \tau)\|_{L^2}^2 \, d\tau - \theta_a(s)^{-(\gamma+1)} \theta_a(s)^{\gamma-1} f(s) \\
    &= \theta_a(s)^{-2} \left[ a \frac{(\gamma + 1)}{2} g_\gamma(s, t) - f(s) \right].
\end{align*}
Thus, although it may not initially be clear how to argue using Gronwall's inequality for a bound on the energy appearing in \eqref{balance}, Neustupa demonstrates in \footnote{Neustupa works with renormalized variables in the non endpoint critical case with $\gamma = 1/3$, without parameters $M$ and $N$.} \cite{Neustupa2014} that one approach is to seek a differential inequality of the form 
\(
\theta_a(s)^2 \partial_s g_\gamma(s, t) + 
A g_\gamma(s, t) \leq  C,
\)
with \( A > 0 \). We will revisit this argument after bounding the terms in the energy balance.

\subsection{Boundedness of the right-hand side terms in the energy balance}
\label{rhsbalance}
We first identify the pressureless terms in the energy balance \eqref{balance}
\begin{align}
\label{np}
    K^I(s, t) =& \int_s^t \frac{\theta_a(\tau)^\gamma}{\theta_a(s)^{\gamma + 1}} \int_{\mathbb{R}^3} \big[ 2|\nabla \Psi_N|^2 |v|^2 + |v|^2 v \cdot \nabla \Psi_N^2 \big] \, dx \, d\tau 
\end{align}
and the pressure term
\begin{equation}
\label{p}
 K^{II}(s, t) = \int_s^t \frac{\theta_a(\tau)^\gamma}{\theta_a(s)^{\gamma + 1}} \int_{\mathbb{R}^3} 2p (v \cdot \nabla \Psi_N^2) \, dx \, d\tau.
\end{equation}
Observe that we can estimate the gradient of $\Psi_N^2$ since
\[
\nabla \Psi_N^2(x, \tau) = \frac{2}{\theta_a(\tau)} \Psi_N(x, \tau) \nabla \varphi_N \left(\frac{x}{\theta_a(\tau)}\right).
\]
Thus, using the notation from the pressure decomposition \eqref{decomposition},
\[
|K^{II}(s, t)| \leq C \int_s^t \frac{\theta_a(\tau)^{\gamma - 1}}{\theta_a(s)^{\gamma + 1}} \int_{\text{supp}(\nabla \Psi_N(\cdot,\tau))} |p v \nabla \Psi_N| \, dx \, d\tau.
\]
\begin{equation}
\label{p2}
    \leq \frac{C^{II}}{\theta_a(s)^{\gamma+1}} \int_s^t \theta_a(\tau)^{\gamma-1} \int\limits_{\text{supp}(\nabla \Psi_N(\cdot,\tau))} \big( |p_1(x, \tau)| + |p_2(x, \tau)| + |p_3(x, \tau)| \big)|v(x, \tau) \Psi_N(x, \tau)| \, dx \, d\tau.
\end{equation}

\paragraph{The \( p_1 \) part.} 
Remember the Lorentz space property
\begin{equation}
\label{normchar}
    \|1_{\text{supp}(\nabla \Psi_N(\cdot,\tau))}\|_{L^{\frac{3}{2},1}(\mathbb{R}^3)} \leq CN^2\theta_a(\tau)^2
\end{equation}
so by Hölder's inequality for Lorentz spaces, using \eqref{suppgradpsi} and the hypothesis \eqref{hip2main},
\[
\int_{\text{supp}(\nabla \Psi_N(\cdot,\tau))} |v(x, \tau)| \, dx \leq CN^2 \theta_a(\tau)^2 \|v(\cdot, \tau)\|_{L^{3,\infty}(\text{supp}(\nabla \Psi_N(\cdot,\tau)))},
\]
\[
\leq CN^2 M \theta_a(\tau)^2.
\]
Then, using this bound and \eqref{kernelp1}, we obtain for all $s \in I_{a,b,c,N}$ and \( t \in (s, 0] \), 
\begin{align*}
\label{p1}
    C^{II} & \int_s^t \theta_a(\tau)^{\gamma-1} \int_{\text{supp}(\nabla \Psi_N(\cdot,\tau))} |p_1(x, \tau)| \, |v(x, \tau)\Psi_N| \, dx \, d\tau \\
    & \leq \frac{C}{N^3} \int_s^t \theta_a(\tau)^{\gamma-4} \left( \int_{ B( (\frac{N}{2}+1) \theta_a(\tau))} |v(x, \tau)|^2 \, dx \right) \left( \int_{\text{supp}(\nabla \Psi_N(\cdot,\tau))} |v(x, \tau)| \, dx \right) d\tau \\
    & \leq \frac{C_\ast M}{N} \int_s^t \theta_a(\tau)^{\gamma-2} \int_{ B( (\frac{N}{2}+1) \theta_a(\tau))} |v(x, \tau)|^2 \, dx \, d\tau,
\end{align*}
where \( C_\ast \in (0, \infty) \) is a universal constant.

\noindent
This small part ($N>>1$) will be absorbed in a similar, yet more optimal, manner compared to \cite{BarkerFernandezPrange2024}, thanks to our choice of \( \eta_N \).

\paragraph{The \( p_2 \) part.} Using the remark \eqref{kernelp2} we obtain for all $s \in I_{a,b,c,N}$ and \( t \in (s, 0] \),
\begin{align*}
    C^{II} & \int_s^t \theta_a(\tau)^{\gamma-1} \int_{\text{supp}(\nabla \Psi_N(\cdot,\tau))} |p_2(x, \tau)| \, |v(x, \tau)\Psi_N| \, dx \, d\tau \\
    & \leq C \int_s^t \theta_a(\tau)^{\gamma-1} \|p_2(\cdot, \tau)\|_{L^{\frac{3}{2},\infty}(\text{supp}(\nabla \Psi_N(\cdot,\tau)))} \|v(\cdot, \tau)\Psi_N\|_{L^{3,1}(\mathbb{R}^3)} \, d\tau \\
    & \leq C\int_s^t \theta_a(\tau)^{\gamma-1} \| v(\cdot, \tau)\|^2_{L^{3,\infty}(B((N+3)\theta_a(\tau)^{1/2}) \setminus B((\frac{N}{2}+1)\theta_a(\tau)))} \|v(\cdot, \tau)\Psi_N\|_{L^{3,1}(\mathbb{R}^3)} \, d\tau.
\end{align*}
Then, by the assumption \eqref{hip2main}, along with the interpolation of \( L^{3,1} \) spaces between \( L^2 \) and \( L^6 \) \cite[Theorem 5.3.1]{BerghLofstrom1976}, the Sobolev inequality, and Young's inequality,
\begin{align*}
    C^{II} & \int_s^t \theta_a(\tau)^{\gamma-1} \int_{\text{supp}(\nabla \Psi_N(\cdot,\tau))} |p_2(x, \tau)| \, |v(x, \tau)\Psi_N| \, dx \, d\tau \\
    & \leq C^{p_2} M^2 \int_s^t \theta_a(\tau)^{\gamma - \frac{3}{2}} \|\Psi_N v(\cdot, \tau)\|_{L^2(\mathbb{R}^3)} \, d\tau + C^{p_2} M^2 \int_s^t \theta_a(\tau)^{\gamma - \frac{1}{2}} \|\nabla(\Psi_N v(\cdot, \tau))\|_{L^2(\mathbb{R}^3)} \, d\tau,
\end{align*}
where \( C^{p_2} \in (0, \infty) \) is a universal constant. 

\noindent
For \( \gamma > 0 \) and \( a \in (0, 4\lambda_S(B_1)) \), we choose \footnote{As we will see, this choice allows us to prove our result with \( a \) spanning the entire interval \( \left(0, 4 \lambda_S(B(1))\right) \).} \( \varepsilon_{a} > 0 \) sufficiently small to ensure that
\begin{equation}
\label{epsilon}
    (2 - \varepsilon_{a}) \lambda_S(B_1) - \frac{a}{2} > \frac{a \varepsilon}{4} > \frac{a\gamma}{2} - \frac{a\gamma}{2 + \varepsilon_{a} / \gamma} > 0.
\end{equation}
Observe that the choice \( \varepsilon_{a} = \frac{4\lambda_S(B_1) - a}{16 \lambda_S(B_1)} \) is feasible because for \( \epsilon > 0 \),
\begin{equation*}
    \frac{a\gamma}{2} - \frac{a\gamma}{2 + \varepsilon / \gamma} =  \frac{a \gamma \varepsilon}{2(2\gamma+ \varepsilon)} < \frac{a \gamma \epsilon}{2(2\gamma)} = \frac{a \varepsilon}{4}
\end{equation*}
and
\begin{equation*}
    (2 - \varepsilon) \lambda_S(B_1) - \frac{a}{2} > \frac{a \epsilon}{4} \iff \varepsilon < 2 \left( \frac{4\lambda_S(B_1) -a}{a+ 4 \lambda_S(B_1)} \right) .
\end{equation*}
Then, using the Cauchy-Schwarz inequality together with Young's inequality, we obtain the following result:
\[
C^{p_2} M^2 \int_s^t \theta_a(\tau)^{\gamma - \frac{3}{2}} \|\Psi_N v(\cdot, \tau)\|_{L^2(\mathbb{R}^3)} \, d\tau
\]
\[
\leq C^{p_2}  M^4 \left(\frac{a\gamma}{2} - \frac{a\gamma}{2 + \varepsilon_{a} / \gamma} \right)^{-1} \frac{\theta_a(s)^{\gamma + 1}}{a(\gamma + 1)} + \frac{1}{2} \left(\frac{a\gamma}{2} - \frac{a\gamma}{2 + \varepsilon_{a} / \gamma} \right) \int_s^t \theta_a(\tau)^{\gamma - 2} \|\Psi_N v(\cdot, \tau)\|_{L^2(\mathbb{R}^3)}^2 d\tau,
\]
where we have used \footnote{This computation is valid for \( \gamma > -1 \).} \( \int_s^0 \theta_a(\tau)^{\gamma-1} =  \theta_a(s)^{\gamma+1} / a(\gamma+1) \), which yields for $\gamma>-1$.

\noindent
Similarly,
\[
C^{p_2} M^2 \int_s^t \theta_a(\tau)^{\gamma - \frac{1}{2}} \|\nabla(\Psi_N v(\cdot, \tau))\|_{L^2(\mathbb{R}^3)} \, d\tau
\]
\[
\leq C_{a, \gamma}  M^4 \theta_a(s)^{\gamma + 1} + \frac{\varepsilon_{a}}{2} \int_s^t \theta_a(\tau)^{\gamma} \|\nabla (\Psi_N v(\cdot, \tau))\|_{L^2(\mathbb{R}^3)}^2 d\tau.
\]

\paragraph{The \( p_3 \) part.} 
Using the observation \eqref{kernelp3} and, once again, the gradient property of \( \eta_N \) from \eqref{suppgradeta}, together with Hölder's inequality for Lorentz spaces and the assumption \eqref{hip1main}, we obtain the following for all \( s \in I_{a,b,c,N} \) and \( t \in (s, 0] \):
\begin{align*}
    C^{II} & \int_s^t \theta_a(\tau)^{\gamma-1} \int_{\text{supp}(\nabla \Psi_N(\cdot,\tau))} |p_3(x, \tau)| \, |v(x, \tau)\Psi_N| \, dx \, d\tau \\
    & \leq CN^2 M \int_s^t \theta_a(\tau)^{\gamma-\frac{1}{2}} \left( \int_{\text{supp} (\nabla \eta_N(\cdot, \tau))} \left( |v|^2 + |p| \right) \, dx \right) d\tau \\
    & \leq C M N^3 \int_s^t \theta_a(\tau)^{\gamma} \left( \int_{\text{supp} (\nabla \eta_N(\cdot, \tau))} \left( |v|^3 + |p|^{\frac{3}{2}} \right) dx \right)^{\frac{2}{3}}.
\end{align*}
Then, recalling the range for \( s \) from \eqref{I}, we obtain the following bound:
\begin{align}
\label{gamma1/2}
    C^{II} & \int_s^t \theta_a(\tau)^{\gamma-1} \int_{\text{supp}(\nabla \Psi_N(\cdot,\tau))} |p_3(\tau)| \, |v(x, \tau)\Psi_N| \, dx \, d\tau \nonumber \\
    & \leq C M N^3 \left( \int_s^t \theta_a(\tau)^{2\gamma} \right)^{\frac{1}{2}}  \left( \int_{-\frac{c}{2}}^0 \left( \int_{B(\frac{b}{4})} \left( |v|^3 + |p|^{\frac{3}{2}} \right) dx \right)^{\frac{4}{3}}  d\tau \right)^{\frac{1}{2}} \nonumber \\
    & \leq  C S M N^3 \frac{ \theta_a(s)^{\gamma+1}  }{ \sqrt{a(\gamma+1)} }  ,
\end{align}
where we have used \footnote{This computation is valid for \( \gamma > -1 \).} \( \left( \int_s^0 \theta_a(\tau)^{2\gamma} \right)^{\frac{1}{2}} =  \theta_a(s)^{\gamma+1} / \sqrt{a(\gamma+1)} \), and the fact that for a suitable Leray-type solution,
\[
S = \left( \int_{-\frac{c}{2}}^0 \left( \int_{B(\frac{b}{4})} \left( |v|^3 + |p|^{\frac{3}{2}} \right) dx \right)^{\frac{4}{3}}  d\tau \right)^{\frac{1}{2}} < \infty.
\]
\noindent
In summary, from the pressure estimates above, we deduce that for all \( s \in I_{a,b,c,N} \) and \( t \in (s, 0] \),
\begin{align}
\label{kii}
    K^{II}(s) \leq& \frac{C_* M}{N} \int_s^t \frac{\theta_a(\tau)^{\gamma-2}}{\theta_a(s)^{\gamma+1}} \|\Psi_N v(\cdot, \tau)\|_{L^2}^2 \, d\tau + \frac{1}{2} \left(\frac{a\gamma}{2} - \frac{a\gamma}{2 + \varepsilon_{a} / \gamma} \right) \int_s^t \frac{\theta_a(\tau)^{\gamma-2}}{\theta_a(s)^{\gamma+1}} \|\Psi_N v(\cdot, \tau)\|_{L^2}^2 \, d\tau \nonumber \\
    &+ \frac{\varepsilon_{a}}{2} \int_s^t \frac{\theta_a(\tau)^{\gamma-2}}{\theta_a(s)^{\gamma+1}} \|\Psi_N v(\cdot, \tau)\|_{L^2}^2 \, d\tau + C_{a,\gamma}(M^4 + SMN^3).
\end{align}
\noindent
For \( K^I(s) \), we observe that
\begin{align*}
    \int_s^t  & \theta_a(\tau)^\gamma \int_{\mathbb{R}^3} |v|^2 v \cdot \nabla \Psi_N^2 \, dx \, d\tau \\
    & \leq C\int_s^t \theta_a(\tau)^{\gamma-1} \| v(\cdot, \tau)\|^2_{L^{3,\infty}(\text{supp}(\nabla \Psi_N(\cdot,\tau)))} \|v(\cdot, \tau) \Psi_N \|_{L^{3,1}(\text{supp}(\nabla \Psi_N(\cdot,\tau))} \, d\tau
\end{align*}
using interpolation and Young inequalities once again, together with our assumptions \eqref{hip1main}  and \eqref{hip2main}, we get 
\begin{align}
\label{ki}
    |K^I(s)| \leq & \frac{1}{2} \left(\frac{a\gamma}{2} - \frac{a\gamma}{2 + \varepsilon_{a} / \gamma} \right) \int_s^t \theta_a(\tau)^{\gamma-2} \|\Psi_N v(\cdot, \tau)\|_{L^2}^2 \, d\tau \nonumber \\
    & + \frac{\varepsilon_{a}}{2} \int_s^t \theta_a(\tau)^{\gamma} \|\nabla(\Psi_N v(\cdot, \tau))\|_{L^2}^2 \, d\tau + C_{a, \gamma} M^4 .
\end{align}

\noindent
Then, from \eqref{balance} combined with \eqref{ki} and \eqref{kii}, we find for \( s \in I_{a,b,c,N} \), and \( t \in (s, 0] \),
\begin{align}
\label{balance2}
    & \frac{\theta_a(t)^\gamma}{\theta_a(s)^{\gamma+1}} \|\Psi_N v(\cdot, t)\|_{L^2}^2 +  \frac{a\gamma}{2 + \varepsilon_{a} / \gamma} \int_s^t \frac{\theta_a(\tau)^{\gamma-2}}{\theta_a(s)^{\gamma+1}} \|\Psi_N v(\cdot, \tau)\|_{L^2}^2 \, d\tau \nonumber \\
    & + (2-\varepsilon_{a} ) \int_s^t \frac{\theta_a(\tau)^{\gamma}}{\theta_a(s)^{\gamma+1}} \|\nabla(\Psi_N v(\cdot, \tau))\|_{L^2}^2 \, d\tau \nonumber \\
    & \leq \frac{1}{\theta_a(s)} \|\Psi_N v(\cdot, s)\|_{L^2}^2 + \frac{C^* M}{N} \int_s^t \frac{\theta_a(\tau)^{\gamma-2}}{\theta_a(s)^{\gamma+1}} \|\Psi_N v(\cdot, \tau)\|_{L^2}^2 \, d\tau + C_{a,\gamma}(M^4 + SMN^3 ),
\end{align}
where the constant \( C^* \in (0, \infty) \) comes from the \( p_1 \) part.

\section{Boundedness of the Scale-Invariant Energy}
\label{boundedness}

The objective of this section is to prove Theorem \ref{premain}.  
\noindent  
We assume\footnote{The proof requires several steps. In particular, it begins by considering \(\gamma > 0\) and proving the boundedness of \( f \). This immediately implies the boundedness of \( g_\gamma \) for \(\gamma > -1\). It is then interesting to observe that the range of \(\gamma\) can be improved in the energy \eqref{balance}, so that \( h_\gamma \) will be bounded for \(\gamma > -1\).} first \(\gamma > 0\).  
We rewrite \eqref{balance2} using the definitions of \( f \) and \( g_\gamma \) provided in \eqref{fg}, as follows:
\begin{align}
\label{balancefg}
    \frac{\theta_a(t)^{\gamma+1}}{\theta_a(s)^{\gamma+1}} f(t) &+  \frac{a\gamma}{2 + \varepsilon_{a} / \gamma} g_\gamma(s, t) + (2 - \varepsilon_{a} ) \int_s^t \frac{\theta_a(\tau)^{\gamma}}{\theta_a(s)^{\gamma+1}} \|\nabla (\Psi_N v)(\cdot, \tau)\|_{L^2}^2 \, d\tau \nonumber\\
    &\leq f(s) + \frac{C^* M}{N} g_\gamma(s, t) + C_{a, \gamma} \big( M^4 + SMN^3 \big).  
\end{align}

\noindent
Then we replace the identity \eqref{fgeq} in \eqref{balancefg} to get for \( s \in I_{a,b,c,N} \), and \( t \in (s, 0] \),
\begin{align}
\label{balancefg2}
    &\frac{\theta_a(t)^{\gamma+1}}{\theta_a(s)^{\gamma+1}} f(t) + \theta_a(s)^2 \frac{\partial g_\gamma}{\partial s}(s, t) + \left( \frac{a \gamma}{2 + \varepsilon_{a} / \gamma} - \frac{a \gamma}{2} - \frac{a}{2} - \frac{C_* M}{N} \right) g_\gamma(s, t) \nonumber \\
    &+ (2-\varepsilon_{a}) \int_s^t \frac{\theta_a(\tau)^{\gamma}}{\theta_a(s)^{\gamma+1}} \|\nabla (\Psi_N v)(\cdot, \tau)\|_{L^2}^2 \, d\tau \nonumber\\
    &\leq  C_{a, \gamma}(M^4 + SMN^3 ).     
\end{align}
We observe that $\left( \frac{a \gamma}{2 + \varepsilon_{a} / \gamma} - \frac{a \gamma}{2} - \frac{a}{2} - \frac{C_* M}{N} \right)$ is negative. Therefore, it is necessary to obtain additional information in damping form to integrate the differential inequality and bound the energy.

\subsection{Gaining Damping with the Gradient Part}
Let us observe that in this subsection we do not use hypothesis \eqref{hip2main}, but only hypothesis \eqref{hip1main}.  
For this part, we follow exactly the computations in~\cite{BarkerFernandezPrange2024}, while verifying control under our assumptions \eqref{hip1main} and \eqref{hip2main}. In this section $s \in I_{a,b,c,N}$ (see \eqref{I}).

\noindent
Let us take \( \xi \in (0, 1) \) and a test function \( \varphi_1^\xi \) defined as:
\[
\varphi_1^\xi(x) = 
\begin{cases} 
1, & \text{if } |x| < 1 + \frac{1}{4}\xi, \\
\in [0, 1], & \text{if } 1 + \frac{1}{4}\xi < |x| < 1 + \frac{3}{4}\xi, \\
0, & \text{if } 1 + \frac{3}{4}\xi < |x|,
\end{cases}
\]
and satisfying the gradient bound
\begin{equation}
\label{gradvarphixi}
    |\nabla \varphi_1^\xi| \leq 4\xi^{-1}.    
\end{equation}
We let \(\varphi_{N,2}^\xi := \varphi_N - \varphi_1^\xi\). Furthermore, we define \(\Psi_1^\xi(x,t) := \varphi_1^\xi\left(\frac{x}{\theta_a(t)}\right)\) and \(\Psi_{N,2}^\xi := \varphi_{N,2}^\xi\left(\frac{x}{\theta_a(t)}\right)\). Thus, it follows that \(\Psi_N = \Psi_1^\xi + \Psi_{N,2}^\xi\).

\noindent
From the control
\begin{align*}
\|\nabla (\Psi_1^\xi v)\|_2^2 
&= \|\Psi_1^\xi \nabla v\|_2^2 + 2 \langle \Psi_1^\xi \nabla v, \nabla (\Psi_1^\xi) \otimes v \rangle_{L^2} + \|\nabla (\Psi_1^\xi) \otimes v\|_2^2 \\
&\leq \|\nabla v\|_{L^2(\text{supp}(\Psi_1^\xi))}^2 + \int_{\mathbb{R}^3} \partial_i \big((\Psi_1^\xi)^2\big) v_j \partial_i v_j \, dx + \|\nabla (\Psi_1^\xi) \otimes v\|_2^2 \\
&\leq \|\nabla (\Psi_N v)\|_2^2 - \frac{1}{2} \int_{\mathbb{R}^3} \Delta(\Psi_1^\xi) |v|^2 \, dx + \|\nabla (\Psi_1^\xi) \otimes v\|_2^2,
\end{align*}
we get
\begin{align}
\label{c1}
    \int_s^t \frac{\theta_a(\tau)^\gamma}{\theta_a(s)^{\gamma+1}} \|\nabla (\Psi_N v)\|_2^2 \, d\tau 
    \geq \int_s^t \frac{\theta_a(\tau)^\gamma}{\theta_a(s)^{\gamma+1}} \|\nabla (\Psi_1^\xi v)\|_2^2 \, d\tau 
    - c_1(s, t, a, \xi),
\end{align}
letting
\begin{align*}
c_1(s, t, a, \xi) := \int_s^t \frac{\theta_a(\tau)^\gamma}{\theta_a(s)^{\gamma+1}} \bigg(
- \frac{1}{2} \int_{\text{supp}(\nabla \Psi_1^\xi)} \Delta (\Psi_1^\xi) |v|^2 \, dx 
+ \|\nabla (\Psi_1^\xi) \otimes v\|_2^2 \bigg) \, d\tau.
\end{align*}
Remembering the formula
\[
\nabla \Psi_1^\xi = \frac{1}{\theta_a(t)} \nabla \varphi_1^\xi \bigg( \frac{x}{\theta_a(t)} \bigg),
\]
we can verify
\[
|c_1(s, t, a, \xi)| \leq \frac{C}{\xi^2} \int_s^t \frac{\theta_a(\tau)^{\gamma-2}}{\theta_a(s)^{\gamma+1}} 
\int_{\text{supp}(\nabla \Psi_1^\xi)} |v(x, \tau)|^2 \, dx \, d\tau,
\]
and then by Hölder's inequality and the assumption \eqref{hip1main}, the constant $c_1$ is well controlled,
\[
-c_1(s, t, a, \xi) \geq -\frac{C}{a \xi^2 (\gamma+1)} M^2.
\]

\subsection{Poincaré inequality with the new localization}
From the fact
\[
\int_{B((1+\xi)\theta_a(\tau))} \nabla (\Psi_1^\xi) \cdot v = 0,
\]
we can apply a right inverse of the divergence operator\footnote{Bogovskii operator in  \cite{Galdi2011} or \cite[Theorem 4]{Borchers1990}.} to obtain the existence of a function
\[
w^\xi(\cdot,\tau) \in W_0^{1,2}(B((1+\xi)\theta_a(\tau))),
\]
for which
\[
\nabla \cdot w^\xi(\cdot,\tau) = \nabla (\Psi_1^\xi) \cdot v(\cdot,\tau)
\]
and
\[
\|\nabla w^\xi\|_{L^2(B((1+\xi)\theta_a(\tau)))} \leq C \|\nabla (\Psi_1^\xi) \cdot v\|_{L^2}.
\]
This \(\Psi_1^\xi v - w^\xi\) is divergence-free in \(B((1+\xi)\theta_a(\tau))\) and have zero trace on the boundary. Then, Poincaré's inequality for divergence-free functions with zero trace, we get
\begin{equation}
\label{tracedivPoincare}
    \|\Psi_1^\xi v - w^\xi\|_{L^2(B((1+\xi)\theta_a(\tau)))} \leq \frac{1}{\sqrt{\lambda_S(B((1+\xi)\theta_a(\tau)))}} \|\nabla (\Psi_1^\xi v - w^\xi)\|_{L^2(B((1+\xi)\theta_a(\tau)))},
\end{equation}
being \(\lambda_S(B(r))\) the first eigenvalue of the Dirichlet-Stokes operator on the ball \(B(r)\). Using homogeneity we also get
\begin{equation}
\label{homoeigenvalue}
    \frac{1}{\sqrt{\lambda_S(B((1+\xi)\theta_a(\tau)))}} = \frac{(1+\xi)\theta_a(\tau)}{\sqrt{\lambda_S(B(1))}}.    
\end{equation}
Thus, using Poincaré's inequality \eqref{tracedivPoincare}, the scaling property \eqref{homoeigenvalue}, and the bound on the test function gradient \eqref{gradvarphixi}, we find
\begin{align*}
\|\Psi_1^\xi v\|_2 
&\leq \|\Psi_1^\xi v - w^\xi\|_{L^2(B((1+\xi)\theta_a(\tau)))} + \|w^\xi\|_2 \\
&\leq \frac{(1 + \xi)\theta_a(\tau)}{\sqrt{\lambda_S(B(1))}} 
\|\nabla (\Psi_1^\xi v - w^\xi)\|_{L^2(B((1+\xi)\theta_a(\tau)))} 
+ \frac{(1 + \xi)\theta_a(\tau)}{\pi} 
\|\nabla w^\xi\|_{L^2(B((1+\xi)\theta_a(\tau)))} \\
&\leq \frac{(1 + \xi)\theta_a(\tau)}{\sqrt{\lambda_S(B(1))}} 
\|\nabla (\Psi_1^\xi v)\|_{L^2(B((1+\xi)\theta_a(\tau)))} 
+ \frac{C(1 + \xi)}{\xi \pi} 
\|v\|_{L^2(B((1+\xi)\theta_a(\tau)) \setminus B(\theta_a(\tau)))}.
\end{align*}
In order to conserve the constant for \(\|\nabla (\Psi_1^\xi v)\|_{L^2(B((1+\xi)\theta_a(\tau)))}\) small, we introduce a small \(\kappa > 0\) such that
\[
\|\Psi_1^\xi v\|_{L^2}^2 \leq \frac{(1+\xi)^2(1+\kappa)\theta_a(\tau)^2}{\lambda_S(B(1))} \|\nabla (\Psi_1^\xi v)\|_{L^2(B((1+\xi)\theta_a(\tau)))}^2 + C_{\xi,\kappa} \|v\|_{L^2(B((1+\xi)\theta_a(\tau)))}^2.
\]
Replacing this control into \eqref{c1} we find 
\begin{equation}
\label{c2}
    \int_s^t \frac{\theta_a(\tau)^\gamma}{\theta_a(s)^{\gamma+1}} \|\nabla (\Psi_N v)\|_{L^2}^2 \, d\tau \geq \frac{\lambda_S(B(1))}{(1+\kappa)(1+\xi)^2} \int_s^t \frac{\theta_a(\tau)^{\gamma-2}}{\theta_a(s)^{\gamma+1}} \|\Psi_1^\xi v\|_{L^2}^2 \, d\tau - c_2(s,t,a,\xi,\kappa),
\end{equation}
letting
\[
c_2(s,t,a,\xi,\kappa) = \frac{C_{\xi,\kappa} \lambda_S(B(1))}{(1+\kappa)(1+\xi)^2} \int_s^t \frac{\theta_a(\tau)^{\gamma-2}}{\theta_a(s)^{\gamma+1}} \|v\|_{L^2(B((1+\xi)\theta_a(\tau)) \setminus B(\theta_a(\tau)) )}^2 \, d\tau + c_1(s,t,a,\xi).
\]
As we have done for \(c_1\), we can control
\[
-c_2(s,t,a,\xi,\kappa) \geq -\frac{C_{\xi,\kappa}}{a(\gamma+1)} M^2.
\]
Then, with the help of the identity
\[
\|\Psi_1^\xi v\|_{L^2}^2 = \|\Psi_N v\|_{L^2}^2 - 2 \langle \Psi_1^\xi v, \Psi_{N,2}^\xi v \rangle_{L^2} + \|\Psi_{N,2}^2 v\|_{L^2}^2
\]
from \eqref{c2} we find
\[
\int_s^t \frac{\theta_a(\tau)^\gamma}{\theta_a(s)^{\gamma+1}} \|\nabla (\Psi_N v)\|_{L^2}^2 \, d\tau \geq \frac{\lambda_S(B(1))}{(1+\kappa)(1+\xi)^2} \int_s^t \frac{\theta_a(\tau)^{\gamma-2}}{\theta_a(s)^{\gamma+1}} \|\Psi_N v\|_{L^2}^2 \, d\tau - c_3(s,t,a,\xi,\kappa).
\]
with
\[
c_3(s, t, a, \xi, \kappa) := \frac{\lambda_S(B(1))}{(1 + \kappa)(1 + \xi)^2} 
\int_s^t \frac{\theta_a(\tau)^{\gamma - 2}}{\theta_a(s)^{\gamma + 1}} 
\left[ 2 (\Psi_1^\xi v, \Psi_{N,2}^\xi v)_2 - \|\Psi_{N,2}^\xi v\|_2^2 \right] d\tau + c_2(s, t, a, \xi, \kappa).
\]

Observing that \(\Psi_{N,2}^\xi(\cdot,t)\) is supported on \(B((N+1)\theta_a(\tau)) \setminus B(\theta_a(\tau))\), from the assumption \footnote{This is the only point where small-order scales (close to the aperture $\sqrt{a}$) are used, that is, we use hypothesis \eqref{hip1main}}\eqref{hip1main} and the Holder's inequality,
\[
-c_3(s, t, \xi, \kappa) \geq -\frac{C_{\xi,\kappa}}{a(\gamma+1)} N M^2,
\]
Thus, we get for \(s \in I_{a,b,c,N}\) and \(t \in (s, 0]\).
\begin{equation}
\label{Poincare}
    \int_s^t \frac{\theta_a(\tau)^\gamma}{\theta_a(s)^{\gamma+1}} \|\nabla (\Psi_N v)(\cdot, \tau)\|_2^2 \, d\tau 
    \geq \frac{\lambda_S(B(1))}{(1 + \kappa)(1 + \xi)^2} g_\gamma(s, t) -\frac{C_{\xi,\kappa}}{a(\gamma+1)} N M^2.
\end{equation}

\subsection{Gronwall's type estimate}
\eqref{gronwall}
The Poincaré-estimate \eqref{Poincare} and the energy control \eqref{balancefg2} implies for \(s \in I_{a,b,c,N}\), and \(t \in (s, 0]\),
\[
\frac{\theta_a(t)^{\gamma+1}}{\theta_a(s)^{\gamma+1}} f(t) + \theta_a(s)^2 \partial_s g_\gamma(s, t) + 
\left( 
\frac{(2 - \varepsilon_{a} ) \lambda_S(B_1)}{(1 + \kappa)(1 + \xi)^2} 
- \frac{a}{2} 
+ \frac{a\gamma}{2 + \varepsilon_{a} / \gamma} 
- \frac{a\gamma}{2} 
- \frac{C_* M}{N} 
\right) g_\gamma(s, t)
\]

\[
\leq C_{a, \gamma} (M^4 + SMN^3 ) + \frac{C_{ \xi, \kappa}}{a(\gamma+1)} N M^2.
\]
We consider \(t = 0\). Defining
\[
A = 
\frac{(2 - \varepsilon_{a} ) \lambda_S(B_1)}{(1 + \kappa)(1 + \xi)^2} 
- \frac{a}{2} 
+ \frac{a\gamma}{2 + \varepsilon_{a} / \gamma} 
- \frac{a\gamma}{2} 
- \frac{C_* M}{N},
\]
and
\[
B := C_{a, \gamma}(M^4 + SMN^3 ) + \frac{C_{ \xi, \kappa}}{a(\gamma+1)} N M^2,
\]
the inequality writes as
\[
\frac{d}{ds} \big(g_\gamma(\cdot, 0)\big)(s) + \frac{A}{\theta_a(s)^2} g_\gamma(s, 0) \leq \frac{B}{\theta_a(s)^2}.
\]
We multiply this expression by the time weight
\begin{equation}
\label{k}
    k(s) = \left(\frac{\theta_a(s)}{\theta_a(s_0)}\right)^{-\frac{2A}{a}},
\end{equation}
for which the integrating factor property fulfills
\[
\frac{dk}{ds} = \frac{A}{\theta_a(s)^2} \left(\frac{\theta_a(s)}{\theta_a(s_0)}\right)^{-\frac{2A}{a}} = \frac{A}{\theta_a(s)^2} k(s),
\]
in order to find
\[
\frac{d}{ds} \big(k g_\gamma(\cdot, 0)\big)(s) \leq \frac{B}{A} \left(\frac{A}{\theta_a(s)^2} k(s)\right).
\]
Integration over \([s_0, s]\) followed by multiplication by \(k^{-1}(s)\) gives
\begin{equation}
\label{gronwall}
    g_\gamma(s, 0) \leq g_\gamma(s_0, 0) \frac{1}{k(s)} + \frac{B}{A} \left(1 - \frac{1}{k(s)}\right).
\end{equation}

\subsection{\texorpdfstring{Boundedness of $g_\gamma$ and $f$}{Boundedness of g_gamma and f}}
\label{boundfg}
From our choice \eqref{epsilon} of $\epsilon_a$, we observe we can take \(  \kappa(a), \xi(a) >0\) small enough
to get 
\[
\frac{(2 - \varepsilon_{a} ) \lambda_S(B_1)}{(1 + \kappa)(1 + \xi)^2} 
- \frac{a}{2} 
+ \frac{a\gamma}{2 + \varepsilon_{a} / \gamma} 
- \frac{a\gamma}{2} >0 
\]
and then the parameter \(N(a,M, \xi(a), \kappa(a))\) large enough to have
\[
A = \frac{(2 - \varepsilon_{a} ) \lambda_S(B_1)}{(1 + \kappa)(1 + \xi)^2} 
- \frac{a}{2} 
+ \frac{a\gamma}{2 + \varepsilon_{a} / \gamma} 
- \frac{a\gamma}{2} 
- \frac{C_* M}{N} > 0.
\]
Then, from the definition of $k$ in \eqref{k} we get, \(1/k(s)\) converges to $0$ when \(s \uparrow 0\), and thus by \eqref{gronwall} we conclude  \(g_\gamma(\cdot, 0)\) is bounded on \([s_0, 0)\).

\noindent
Since \(f\) is independent of \(\gamma\), the boundedness of \(f(s)\) can be established by setting \(\gamma = 1\) and relying on the control provided by \(g_1\). Consider any \(t_1 \in [s_0/2, 0)\) and define \(s_1 := 2t_1\). This ensures \(s_1 < t_1 < 0\) and satisfies the relation 
\[
2(t_1 - s_1) = -s_1 = \frac{\theta_a(s_1)^2}{a}.
\]
Additionally, for \(\tau \in (s_1, t_1)\), the inequality
\[
\frac{1}{\sqrt{2}} \theta_a(s_1) = \theta_a(t_1) < \theta_a(\tau) < \theta_a(s_1) = \sqrt{2} \theta_a(t_1),
\]
holds. Consequently, we obtain:
\[
\frac{1}{t_1 - s_1} \int_{s_1}^{t_1} f(\tau) d\tau \leq \frac{2a}{\theta_a(s_1)^2} \int_{s_1}^{t_1} f(\tau) d\tau \leq 2a \int_{s_1}^{0} \frac{1}{\theta_a(s_1)^2} f(\tau) d\tau 
= 2a g_1(s_1, 0).
\]
Moreover, there exists \(s_1' \in (s_1, t_1)\) such that
\begin{equation}
\label{auxf}
    f(s_1') \leq \frac{1}{t_1 - s_1} \int_{s_1}^{t_1} f(\tau) d\tau \leq 2a \sup_{s \in [s_0, 0)} g_1(s, 0).  
\end{equation}
By combining inequality \eqref{balancefg}, evaluated at \(t = t_1\) and \(s = s_1'\), with \eqref{auxf}, we obtain:
\begin{align*}
\frac{1}{2} f(t_1) 
&\leq \frac{\theta_a(t_1)^2}{\theta_a(s_1')^2} f(t_1) \\
&\leq f(s_1') + \frac{C_* M}{N} g_1(s_1', 0) 
+ C_{a, \gamma} (M^4 + SMN^3 ) \\
&\leq 2a \sup_{s \in [s_0, 0]} g_1(s, 0) 
+ \frac{C_* M}{N} g_1(s_1', 0) 
+ C_{a, \gamma} (M^4 + SMN^3 ).
\end{align*}
Hence, \(f\) remains bounded on a small, non-empty interval \([s_0/2, 0)\), and subsequently on \((-c, 0)\) by leveraging the boundedness of the energy for \(s \in (-c, s_0/2)\). Consequently, we conclude:
\[
\esssup_{s \in (-1, 0)} \frac{1}{\theta_a(s)} \|\Psi_N v(\cdot, s)\|_2^2 < +\infty,
\]
which implies the boundedness of $g_\lambda$ for $\lambda>-1$,
\[
\esssup_{s \in (-c,0)}  \int_s^0 \frac{\theta_a(\tau)^{\gamma-2}}{\theta_a(s)^{\gamma+1}} \|\Psi_N v(\tau)\|_2^2 d\tau < +\infty.
\]

\subsection{\texorpdfstring{Boundedness of $h_\gamma$}{Boundedness of h_gamma}}
Utilizing the boundedness of \(f\) and \(g_\gamma\), we can now revisit the computations that yield the estimates for \(K^{I}\) and \(K^{II}\), this time without the need to absorb terms, to obtain
\begin{align*}
    |K^{II}(s)| \leq& \frac{C_* M}{N} \int_s^t \frac{\theta_a(\tau)^{\gamma-2}}{\theta_a(s)^{\gamma+1}} \|\Psi_N v(\cdot, \tau)\|_{L^2}^2 \, d\tau + C \int_s^t \frac{\theta_a(\tau)^{\gamma-2}}{\theta_a(s)^{\gamma+1}} \|\Psi_N v(\cdot, \tau)\|_{L^2}^2 \, d\tau \nonumber \\
    &+ \frac{1}{2} \int_s^t \frac{\theta_a(\tau)^{\gamma-2}}{\theta_a(s)^{\gamma+1}} \|\Psi_N v(\cdot, \tau)\|_{L^2}^2 \, d\tau + C_{a, \gamma} \left( M^4 + SMN^3 \right)
\end{align*}
and
\begin{align*}
    |K^I(s)| \leq & C \int_s^t \theta_a(\tau)^{\gamma-2} \|\Psi_N v(\cdot, \tau)\|_{L^2}^2 \, d\tau \nonumber \\
    & + \frac{1}{2} \int_s^t \theta_a(\tau)^{\gamma} \|\nabla(\Psi_N v(\cdot, \tau))\|_{L^2}^2 \, d\tau + C_{a,\gamma}M^4 .
\end{align*}
These computations\footnote{Since \(f\) is bounded, there is no longer a need to select \(\varepsilon_a\) or assume \(\gamma > 0\).} are valid for \(\gamma > -1\), and from the balance \eqref{balance}, we obtain:
\begin{align*}
    \esssup_{s \in I_{a,b,c,N}} & \frac{1}{\theta_a(s)} \int_s^0 \|\nabla(\Psi_N v)\|_2^2 d\tau \\
    \leq & \esssup_{s \in I_{a,b,c,N}} \left(\frac{1}{\theta_a(s)} \|\Psi_N v(s)\|_2^2 + \left( C + \frac{C^* M}{N} - \frac{a\gamma}{2} \right) \int_s^0 \frac{\theta_a(\tau)^{\gamma-2}}{\theta_a(s)^{\gamma+1}} \|\Psi_N v(\tau)\|_2^2 d\tau \right) \\
    & + C \left( M^4 + SMN^3 \right).
\end{align*}
which implies $h_\gamma$ is bounded for $\gamma>-1$ and we conclude the proof of Theorem \ref{premain}.

\subsection{Proof of Part B) in Theorem \ref{premain}}
We need to recalibrate our dynamic decomposition. Now, we consider
\begin{equation*}
    \eta_N(x,\tau) = \varphi_{2N+1} \left( \frac{x}{\theta_a(\tau)} \right)
\end{equation*}
Then, we have 
\begin{align*}
\label{suppgradetapri}
    \text{supp} (\nabla \eta_N(\cdot, \tau)) &\subset B \left( (2N+1) \theta_a(\tau) \right) \setminus B \left( (2N+2) \theta_a(\tau)\right)
\end{align*}
with the following bounds: \( |\nabla \eta_N| \lesssim \frac{C}{N \theta(t)} \) and \( |\nabla^2 \eta_N| \lesssim \frac{C}{N^2 \theta(t)^2} \). 

\noindent
For the time variable $s$ in the set
\begin{equation}
\label{Ipri}
    s \in  (-\frac{c}{2},0) \cap (-\frac{1}{2a},0) \cap \left(-\frac{b^2}{a 4^2 (2N+2)^2}, 0\right) =: I_{a,b,c,N},
\end{equation}
which implies for $\tau \in (s,0)$,
\begin{equation*}
    B((N+1)\theta_a(\tau)) \subset B((2N+2)\theta_a(\tau)) \subset B(\frac{b}{4}),
\end{equation*}
we focus on the boundedness of
\[
K^{II}(s, t) \int_s^t \frac{\theta_a(\tau)^\gamma}{\theta_a(s)^{\gamma + 1}} \int_{\mathbb{R}^3} 2p (v \cdot \nabla \Psi_N^2) \, dx \, d\tau
\]
\begin{equation*}
\label{p2cri}
    \leq \frac{C^{II}}{\theta_a(s)^{\gamma+1}} \int_s^t \theta_a(\tau)^{\gamma-1} \int\limits_{\text{supp}(\nabla \Psi_N(\cdot,\tau))} \big( | p_1(x, \tau)| + |p_2(x, \tau)| + |p_3(x, \tau)| \big)|v\Psi_N| \, dx \, d\tau.
\end{equation*}
where
\[
p_1(x, \tau) := \frac{1}{4\pi} \int_{B((\frac{N}{2}+1))\theta_a(\tau))} \frac{\partial^2}{\partial y_i \partial y_j} \left( \frac{1}{|x-y|} \right) [\eta v_i v_j](y, \tau) \, dy,
\]
\[
p_2(x, \tau) := \frac{1}{4\pi} \int_{B\left((2N+2) \theta_a(\tau)\right) \setminus B((\frac{N}{2}+1))\theta_a(\tau))} \frac{\partial^2}{\partial y_i \partial y_j} \left( \frac{1}{|x-y|} \right) [\eta v_i v_j](y, \tau) \, dy,
\]
and
\begin{align*}
    p_3(x, \tau) :=& \frac{1}{2\pi} \int_{\text{supp} (\nabla \eta_N(\cdot, \tau))} \frac{x_i - y_i}{|x-y|^3} \left( \frac{\partial \eta}{\partial y_j} v_i v_j \right)(y, \tau) \, dy \\
    & + \frac{1}{4\pi} \int_{\text{supp} (\nabla \eta_N(\cdot, \tau))} \frac{1}{|x-y|} \left( \frac{\partial^2 \eta}{\partial y_i \partial y_j} v_i v_j \right)(y, \tau) \, dy \\
    & + \frac{1}{2\pi} \int_{\text{supp} (\nabla \eta_N(\cdot, \tau))} \frac{x_i - y_i}{|x-y|^3} \left( \frac{\partial \eta}{\partial y_j} p \right)(y, \tau) \, dy \\
    &+ \frac{1}{4\pi} \int_{\text{supp} (\nabla \eta_N(\cdot, \tau))} \frac{1}{|x-y|} (\Delta \eta p)(y, \tau) \, dy.
\end{align*}
Observe that the bound for the term \( p_1 \) remains the same:
\begin{align*}
\label{p1cri}
    C^{II} & \int_s^t \theta_a(\tau)^{\gamma-1} \int_{\text{supp}(\nabla \Psi_N(\cdot,\tau))} |p_1(x, \tau)| \, |v(x, \tau)\Psi_N| \, dx \, d\tau \\
    &\leq \frac{C_\ast M}{N} \int_s^t \theta_a(\tau)^{\gamma-2} \int_{ B( (\frac{N}{2}+1) \theta_a(\tau))} |v(x, \tau)|^2 \, dx \, d\tau.
\end{align*}

\noindent
For \( p_2 \), by Calderón-Zygmund's theorem and \eqref{hip2maincri},
\begin{equation*}
\label{kernelp2cri}
    \|p_2(\cdot, \tau)\|_{L^{\frac{3}{2},\infty}(\text{supp}(\nabla \Psi_N(\cdot,\tau)))} \leq C \| v(\cdot, \tau)\|^2_{L^{3,\infty}(B((2N+2)\theta_a(\tau)) \setminus B((\frac{N}{2}+1)\theta_a(\tau)))} \leq CM^2.
\end{equation*}
Then, by choosing \( \varepsilon_a \) as in \eqref{epsilon} and proceeding as before under our new hypothesis \eqref{hip1maincri} and \eqref{hip2maincri}, we obtain a similar bound:
\begin{align*}
    C^{II} & \int_s^t \theta_a(\tau)^{\gamma-1} \int_{\text{supp}(\nabla \Psi_N(\cdot,\tau))} |p_2(x, \tau)| \, |v(x, \tau)\Psi_N| \, dx \, d\tau \\
    & \leq \frac{1}{2} \left(\frac{a\gamma}{2} - \frac{a\gamma}{2 + \varepsilon_{a} / \gamma} \right) \int_s^t \frac{\theta_a(\tau)^{\gamma-2}}{\theta_a(s)^{\gamma+1}} \|\Psi_N v(\cdot, \tau)\|_{L^2}^2 \, d\tau \nonumber \\
    &+ \frac{\varepsilon_{a}}{2} \int_s^t \frac{\theta_a(\tau)^{\gamma-2}}{\theta_a(s)^{\gamma+1}} \|\Psi_N v(\cdot, \tau)\|_{L^2}^2 \, d\tau + C_{a,\gamma} M^4.
\end{align*}
For \( p_3 \), we observe that 
\begin{align*}
    \text{dist} & \big( \text{supp}(\nabla \Psi_N(\cdot,\tau)), \text{supp} (\nabla \eta_N(\cdot, \tau)) \big) \nonumber\\
    &\geq (2N+1) \theta_a(\tau)- (N+1) \theta_a(\tau) \nonumber\\
    &\geq N \theta_a(\tau),
\end{align*}
which implies
\begin{equation*}
\label{kernelp3cri}
    |p_3(x, \tau)| \leq \frac{C}{N^3 \theta(\tau)^3} \int_{\text{supp} (\nabla \eta_N(\cdot, \tau))} \big(|v|^2 + |p| \big)(y, \tau) \, dy.
\end{equation*}
Hence, by utilizing the weak critical hypothesis \eqref{weakcrip} on the pressure \( p \), we obtain the following for all \( s \in I_{a,b,c,N} \) and \( t \in (s, 0] \):
\begin{align*}
    C^{II} & \int_s^t \theta_a(\tau)^{\gamma-1} \int_{\text{supp}(\nabla \Psi_N(\cdot,\tau))} |p_3(x, \tau)| \, |v(x, \tau)\Psi_N| \, dx \, d\tau \\
    & \leq \frac{C M}{N} \int_s^t \theta_a(\tau)^{\gamma-2} \left( \int_{\text{supp} (\nabla \eta_N(\cdot, \tau))} \left( |v|^2 + |p| \right) \, dx \right) d\tau \\
    & \leq C M^2 \int_s^t \theta_a(\tau)^{\gamma-1} \\
    & = C M^2  \frac{\theta_a(s)^{\gamma+1}} {a(\gamma+1)}.
\end{align*}
Thus, the right-hand side of the energy balance \eqref{balance} can be controlled for \( \gamma > 0 \) as follows:
\begin{align*}
\label{balancefgcri}
    \frac{\theta_a(t)^{\gamma+1}}{\theta_a(s)^{\gamma+1}} f(t) &+  \frac{a\gamma}{2 + \varepsilon_{a} / \gamma} g_\gamma(s, t) + (2 - \varepsilon_{a} ) \int_s^t \frac{\theta_a(\tau)^{\gamma}}{\theta_a(s)^{\gamma+1}} \|\nabla (\Psi_N v)(\cdot, \tau)\|_{L^2}^2 \, d\tau \nonumber\\
    &\leq f(s) + \frac{C^* M}{N} g_\gamma(s, t) + C_{a, \gamma}  M^4 \big.  
\end{align*}
It is interesting to observe that our control in the right hand side does not grow as \( N \to \infty \).

\noindent
Now, once again replacing \( f(s) \) in terms of \( g_\lambda \) as indicated in \eqref{fgeq}, and applying Poincaré's inequality, we obtain:
\[
\frac{\theta_a(t)^{\gamma+1}}{\theta_a(s)^{\gamma+1}} f(t) + \theta_a(s)^2 \partial_s g_\gamma(s, t) + 
\left( 
\frac{(2 - \varepsilon_{a} ) \lambda_S(B_1)}{(1 + \kappa)(1 + \xi)^2} 
- \frac{a}{2} 
+ \frac{a\gamma}{2 + \varepsilon_{a} / \gamma} 
- \frac{a\gamma}{2} 
- \frac{C_* M}{N} 
\right) g_\gamma(s, t)
\]
\[
\leq C_{a, \gamma} M^4  + \frac{C_{ \xi, \kappa}}{a(\gamma+1)} N M^2,
\]
where \( \xi \) and \( \kappa \) can be taken arbitrarily small. Proceeding as in Section \ref{gronwall} and Section \ref{boundfg}, for \( k \) defined in \eqref{k} and \( s_0 \) chosen in \eqref{Ipri}, we obtain:
\begin{equation*}
\label{gronwallcri}
    g_\gamma(s, 0) \leq g_\gamma(s_0, 0) \frac{1}{k(s)} + \frac{B}{A} \left(1 - \frac{1}{k(s)}\right),
\end{equation*}
with
\[
A = 
\frac{(2 - \varepsilon_{a} ) \lambda_S(B_1)}{(1 + \kappa)(1 + \xi)^2} 
- \frac{a}{2} 
+ \frac{a\gamma}{2 + \varepsilon_{a} / \gamma} 
- \frac{a\gamma}{2} 
- \frac{C_* M}{N} > 0
\]
and
\[
B := C_{a, \gamma} M^4 + \frac{C_{ \xi, \kappa}}{a(\gamma+1)} N M^2,
\]

Thus, the function \( \frac{1}{N} g_\gamma \) is bounded by a constant that depends on the parameters \( M \) and \( a \), but not on \( N \), over the interval \( I_{a,b,c,N} \). Following the computations used to bound the \( N \)-dependent functions \( f \) and \( h_\gamma \), we see that similar bounds hold for \( \frac{1}{N} f \) and \( \frac{1}{N} h_\gamma \) over the interval \( I_{a,b,c,N} \). This completes the proof of Part B).

\section{Proof of Corollary \ref{inter}
 and \ref{main}}
\label{conclusions}

\subsection{Corollary \ref{inter}}
The fact that the boundedness of the function \( f \) implies the boundedness of the functions \( g_\gamma \) and \( h_\gamma \) follows directly as a corollary of the proof of Theorem \ref{premain}. It remains to establish the boundedness of the \( L^3 L^3 \) critical quantity.
By interpolation between \(L^2\) and \(L^6\), we obtain
\begin{align*}
    &\frac{ 1 }{ - s } \int_{s}^{0} \int_{ \text{supp} (\Psi_N)} | \Psi_N {\bf v}  |^3(x,\tau)  \, dx \, d \tau, \\
    \leq& \frac{a}{\theta(s)^2} \int_{s}^{t_0} \left(  \theta (\tau)^{-\frac{1}{4}}   \| \Psi_N {\bf v} \|_2 ^{\frac{3}{2}}    \right) \left(  \theta (\tau)^{\frac{1}{4}} \| \Psi_N {\bf v} \|_6 ^{\frac{3}{2}}   \right)  \, d\tau \nonumber \\
    \leq& \frac{2a}{3^{\frac{3}  {4}} \pi \theta (s)^2 }   \left( \int_{s}^{t_0}   \theta (\tau)^{-1}   \| \Psi_N {\bf v} \|_2 ^{6}  \, d\tau \right)^\frac{1}{4} \left( \int_{s}^{t_0}  \theta (\tau )^\frac{1}{3}  \| \nabla (\Psi_N {\bf v}) \|_2 ^{2}  \, d \tau \right)^\frac{3}{4}, \nonumber
\end{align*}
so that
\begin{align*}
    &\frac{ 1 }{ - s } \int_{s}^{0} \int_{ \text{supp} (\Psi_N)} | \Psi_N {\bf v}  |^3(x,\tau)  \, dx \, d \tau, \\
    \leq& \frac{2 a }{3^{\frac{3}  {4}} \pi }   \left(  \left[ \underset{ s < \tau < 0 }{\mathrm{ess \,  sup}} \frac{ \theta ( \tau ) ^ {\frac{1}{3} } } {\theta(s) ^{\frac{4}{3} } } \| \Psi_N {\bf v} ( \cdot , \tau ) \|_2 ^{2} \right]^2 \, \int_{s}^{0} \frac{ \theta  (\tau )^{-\frac{5}{3}}}{ \theta(s) ^{\frac{4}{3}}   } \| \Psi_N {\bf v} \|_2 ^{2}  \, d \tau   \right)^\frac{1}{4} \nonumber \\
    & \times \left( \int_{s}^{0}  \frac{ \theta ( \tau)^{\frac{1}{3}}}{ \theta(s)^{\frac{4}{3}}} \| \nabla  ( \Psi_N {\bf v} ) \|_2^2  \, d \tau  \right)^\frac{3}{4}.
\end{align*}
Hence, employing the boundedness of $f$, $g_\gamma$ and $h_\gamma$ for $\gamma = \frac{1}{3}$, we get
\begin{equation}
\label{l3l3}
    \esssup_{s \in I_{a,b,c,N}} \frac{1}{\theta_a(s)^2} \int_s^0 \int |\Psi_N v(x,\tau)|^3 dx d\tau < C_N < +\infty.
\end{equation}
Then, by the pigeonholing principle we obtain a sequence of times \(t_k \uparrow 0 \) such that $\int_{B(N \theta_a(t_k))} | v(\cdot, t_k) |^3 $ is bounded and Corollary \ref{inter} is proved.

\subsection{Corollary \ref{main}}
Utilizing \eqref{l3l3} and \eqref{hipgen}, we have
\[
\esssup_{s \in (-c,0))} \frac{1}{\theta_a(s)^2} \int_s^0 \int_{B(b)} |v(x, \tau)|^3 dx d\tau < +\infty.
\]
Then, the pigeonholing argument lead to find a sequence of times \(t_k \uparrow 0\) for which
\[
\sup_k \int_{B(b)} |v(x, t_k)|^3 dx < +\infty.
\]
By \cite[Theorem 1.1]{AlbrittonBarker2020} we get that \((0,0)\) is a regular point.

\appendix

\section{\texorpdfstring{Scaling invariance of $f$, $g_\gamma$ and $h_\gamma$}{Scaling invariance of f, g_gamma and h_gamma}}
\label{critical}
In this section, we specify the scaling invariance property for the principal functions under consideration. Although \(f\), \(g_\gamma\), and \(h_\gamma\) inherently depend on \(N\), we omit this subscript as the scaling invariance is intrinsic to the paraboloid and remains independent of \(N\).  
Denoting:
\begin{equation*}
    g_v(s,t) = \int_s^t\frac{\theta_a(\tau)^{\gamma-2}}{\theta_a(s)^{\gamma+1}}\|\Psi_N v (\cdot,s)\|_{L^2}^2d\tau
\end{equation*}
and 
\begin{equation*}
    v_\lambda (x, \tau) = \lambda^{-1} v (\lambda^{-1}x, \lambda^{-2}\tau),
\end{equation*}
we obtain by performing a change of variables
\begin{align*}
    g_{v_\lambda}(s,t) &= \int_s^t \frac{\theta_a(\tau)^{\gamma-1}}{\theta_a(s)^{\gamma+1}} \| \Psi_N v_\lambda(\cdot,\tau) \|_{L^2}^2 d\tau \\   
    & = \int_{ s / \lambda^2 }^{t / \lambda^2} \frac{\theta_a(\lambda^2\tau)^{\gamma-2}}{\theta_a(s)^{\gamma+1}} \int_{\mathbb{R}^3} \left| \varphi_{N} \left( \frac{\lambda x}{\sqrt{-a \lambda^2 \tau}} \right) \frac{1}{\lambda} v ( x, \tau) \right|^2 \, \lambda^3 dx \, \lambda^2 d\tau \\
    & = \int_{s / \lambda^2}^{t / \lambda^2} \frac{\theta_a(\tau)^{\gamma-2}}{\theta_a (s / \lambda^2)^{\gamma+1}} \|\Psi_{N} v (\cdot,\tau)\|_{L^2}^2 d\tau \\
    & = g_{v}(s / \lambda^2,t / \lambda^2).
\end{align*}
Similarly, for $h_\gamma$, since $ \nabla ( \Psi_N v_\lambda ) ( \lambda x , \lambda^2 \tau )  = \frac{1}{\lambda} \nabla ( \Psi_N v_\lambda (\lambda \cdot , \lambda^2 \tau)  ) (\lambda x) $, applying the same change of variables in the integral, we obtain:
\begin{align*}
    h_{v_\lambda}(s,t) &:= \int_s^t \frac{\theta_a(\tau)^{\gamma}}{\theta_a(s)^{\gamma+1}}  \| \nabla ( \Psi_N v_\lambda ) (\cdot,\tau)  \|_{L^2}^2 d\tau \\  
    & = \int_{ s / \lambda^2 }^{t / \lambda^2} \frac{\theta_a(\lambda^2\tau)^{\gamma}}{\theta_a(s)^{\gamma+1}} \int_{\mathbb{R}^3} \left| \frac{1}{\lambda} \nabla \left( \varphi_{N} \left( \frac{\lambda \cdot }{\sqrt{-a \lambda^2 \tau}} \right) \frac{1}{\lambda} v ( \cdot, \tau) \right) \right|^2 \, \lambda^3 dx \, \lambda^2 d\tau \\
    & = \int_{s / \lambda^2}^{t / \lambda^2} \frac{\theta_a(\tau)^{\gamma}}{\theta_a (s / \lambda^2)^{\gamma+1}} \| \nabla ( \Psi_{N} v (\cdot,\tau) ) \|_{L^2}^2 d\tau \\
    & = g_{v}(s / \lambda^2,t / \lambda^2)
\end{align*}
and
\begin{equation*}
    f_{v_\lambda}(s) := \frac{1}{\theta_a(s)} \int \left| \Psi_N v_\lambda (\cdot,s) \right|^2 dx = \frac{1}{\theta_a(s)} \int \left| \varphi_{N} \left( \frac{\lambda x}{\sqrt{-a \lambda^2 \tau}} \right) \frac{1}{\lambda} v (\cdot,s / \lambda^2) \right|^2 \lambda^3 dx = f_{v}(s / \lambda^2).
\end{equation*}

\section*{Conflict of Interest}
The author declare that he have no conflict of interest.

\section*{Data Availability Statement}
Data sharing is not applicable to this article as no datasets were generated or analyzed during the current study.

\section*{Acknowledgment}
PF is supported by the Basque Government through the BERC 2022-2025 program and by the Spanish State Research Agency through BCAM Severo Ochoa CEX2021-001142

\end{document}